\documentclass[journal]{IEEEtran}
\usepackage{amssymb}
\usepackage{epsfig}
\usepackage{graphicx}
\usepackage{colortbl}

\newtheorem{theorem}{Theorem}

\newtheorem{remark}{Remark}
\newtheorem{lemma}{Lemma}

\newtheorem{definition}{Definition}

\newcommand{\hs}{\hspace}

\ifCLASSINFOpdf
\else
\fi
\begin{document}
%
\title{Consensus Control of Multi-agent Systems with Optimal Performance}
%
%
%

\author{Juanjuan~Xu,~Huanshui~Zhang~
\thanks{*This work is supported by the National Natural
Science Foundation of China under Grants 61573221, 61633014. H. Zhang is the corresponding author.}
\thanks{J. Xu is with School of Control Science and Engineering, Shandong University, Jinan, Shandong, P.R. China 250061.
        {\tt\small juanjuanxu@sdu.edu.cn}}
\thanks{H. Zhang is with School of Control
Science and Engineering, Shandong University, Jinan, Shandong, P.R.China 250061.
        {\tt\small hszhang@sdu.edu.cn}}
}

\maketitle

\begin{abstract}
The consensus control with optimal cost remains major challenging  although consensus control problems have been well studied in recent years.
In this paper, we study the consensus control of multi-agent system associated with a given cost function.
The main contribution is to present the distributed control protocol while minimizing the given cost function with positive semi-definite weighting matrix of control.
The derived controller is composed of two parts: One part is the feedback of the individual state which minimizes the cost function and the
other part is the feedback of the relative states between the neighbours which guarantees the consensus. The presented results are new to the best of our knowledge.

\end{abstract}


%
\IEEEpeerreviewmaketitle

\section{Introduction}

Consensus of multi-agent systems have wide applications in plenty of areas such as automata theory
\cite{liao}, distributed computation \cite{lynch}, formation control \cite{fax}
and so on. \cite{olfati}-\cite{hatano} studied the consensus problem of first-order integrators under various
communication topologies, including directed or undirected, fixed or switching, and random information network.
The mean-square consensus for multi-agent systems consisting of second-order integrator was investigated in \cite{tian}.
\cite{linlisun} further studied the consensus problems of a class of high-order
multi-agent systems. \cite{ma} and \cite{you} obtained the necessary and sufficient condition for the consensusability of
the discrete-time and continuous-time multi-agent systems respectively.

More recently, optimal consensus has become one of the most hot topics and attracted
increasing attention in various scientific communities. There have been much research which were aimed to find
the optimal protocols under a certain pre-specified cost function \cite{cao}-\cite{wanglong}. The fundamental problem in optimal
consensus problem is that the global information is generally needed to achieve the global optimization, which is difficult due to
the communication constraints and the limitations of sensor. Taking account of the existing work, \cite{cao}-\cite{wanglong}
derived the optimal consensus protocol by solving a Riccati equation. However, the distributed structure of the optimal controller stresses
strict conditions for the communication topology. To be specific, in \cite{cao} and \cite{fangfangzhang}, the graph is required to be complete.
In \cite{wanglong}, the leader-follower topology is restricted to be a star graph, which
requires that each agent has a directed path to the leader. \cite{kristian}-\cite{zhangfeng} studied the inverse optimality problem
of multi-agent systems. In the inverse optimality problem, the cost function and the distributed optimal protocol
are needed to be designed simultaneously. In addition, the design of the distributed protocol to optimize local cost function has also been studied in
\cite{bauso}-\cite{semsar}. \cite{ngauyen} and \cite{gupta} proposed two consensus protocols by considering the linear
quadratic regulator. However, the consensus controller is suboptimal.


For most work in the literature, the weighting matrix of control in the cost functions is constrained to be positive definite.
This leads to the unique existence of the optimal controller to minimize the cost function and also makes a restriction for the
structure of the optimal controller which uses all agents' information in general. Thus, the simultaneous controller design of
the distribution and optimization is a difficult issue for the case with positive definite weighting matrix of control.


In this paper, we will consider the consensus control of multi-agent system with optimization performance. Since there is seldom optimal and
distributed control protocol for the cost function with a positive definite weighting matrix of control,
we will study a class of cost functions with positive
semi-definite weighting matrix of control.
The optimal control protocol is given to minimize the given cost function which also has a distributed structure.
To be specific, the derived controller is consisting of two parts where one part is the feedback of the individual state which minimizes the cost function and the other part is the feedback of the relative states between the neighbours which guarantees the consensus.

The remainder of the paper is organized as follows. Section II presents the problem studied in this paper.
Section III shows some preliminary knowledge about algebraic graph theory and the result of the optimization problem for single-agent
system. The optimal consensus controller is given in Section IV. Numerical examples are given in Section V. Some concluding words are given in
the last section.

The following notation will be used throughout this paper: $R^n$
denotes the family of $n$-dimensional vectors; $x'$ denotes the
transpose of $x$; a symmetric matrix $M>0~ (\geq 0)$ means that $M$ is
strictly positive-definite (positive semi-definite).


\section{Problem Formulation}

Consider the multi-agent system
\begin{eqnarray}
\dot{x}_i(t)&=&Ax_i(t)+Bu_i(t),i=1,\cdots,N,\label{r1}
\end{eqnarray}
where $x_i\in R^n$ is the state of the $i$th agent, $u_i\in R^m$ is the control input of the $i$th agent. $A,B$ are constant
matrices with compatible dimensions. The initial values are given by $x_i(0), i=1,\cdots,N.$
The cost function for agent $i$ is given by
\begin{eqnarray}
J_i&=&\int_0^\infty\Big(x_i'(t)Qx_i(t)+u_i'(t)Ru_i(t)\Big)dt,\label{r2}
\end{eqnarray}
where $Q\geq0,R\geq0.$

\begin{definition}
The multi-agent system (\ref{r1}) with a fixed and connected undirected graph is said to be consensusable if for any finite initial values $x_i(0),$
there exists a distributed control protocol
such that the multi-agent achieves consensus, i.e.
\begin{eqnarray}
\lim_{t\rightarrow\infty}E\|x_j(t)-x_i(t)\|^2=0,~\forall i, j=1,\cdots,N.\nonumber
\end{eqnarray}
\end{definition}

The aim of the paper is to find a distributed control protocol which guarantees the consensus of the multi-agent system (\ref{r1}) and
minimizes the cost function $\sum_{i=1}^NJ_i$.

\begin{remark}
In the literature, the weighting matrix $R$ is assumed to be positive definite. Thus the optimal
controller exists uniquely. However, the optimal controller may not always have a distributed structure.
Thus the simultaneous design of optimization and distribution remains challenging. In this paper, by considering the
optimization problem with a positive semi-definite weighting matrix of control, we will propose a consensus
protocol with optimization performance.

\end{remark}

Noting the semi-positive definiteness of the weighting matrices, we will apply the Moore-Penrose inverse of a matrix. To this end, we firstly present some definitions of the Moore-Penrose inverse. From \cite{pinv}, for
a given matrix $M\in R^{n\times m}$, there exists a unique matrix in
$R^{m\times n}$ denoted by $M^{\dag}$ such that
\begin{eqnarray}
MM^{\dag}M=M,~M^{\dag}MM^{\dag}=M^{\dag},\nonumber\\
(MM^{\dag})¡¯=MM^{\dag},~(M^{\dag}M)'=M^{\dag}M.\nonumber
\end{eqnarray}
The matrix $M^{\dag}$ is called the Moore-Penrose inverse of $M$. The following lemma is from \cite{rami}.
\begin{lemma}\label{lem}
Let matrices $L, M$ and $N$ be given with appropriate size. Then, $LXM=N$ has a solution $X$ if and only if $LL^{\dag}NMM^{\dag}=N$. Moreover, the solution of $LXM=N$ can be expressed as
$X=L^{\dag}NM^{\dag}+Y-L^{\dag}LY MM^{\dag},$ where $Y$ is a matrix with appropriate size.
\end{lemma}
In particular, let $M=I$, we have $LX=N$ has a solution if and only
if $LL^{\dag}N = N$. This is also equivalent to $Range(N)\subseteq
Range(L)$ where $Range(N)$ is the range of $N.$


\section{Preliminaries}

\subsection{Algebraic Graph Theory}
In this paper, the information exchange among agents is modeled by an undirected graph.
Let $\mathcal{G}=(\mathcal{V}, \mathcal{E}, \mathcal{A})$ be a diagraph with the set of vertices
$\mathcal{E}=\{1, \ldots, N\}$, the set of edges $\mathcal{E}\subset \mathcal{V}\times \mathcal{V}$,
and the weighted adjacency matrix $\mathcal{A}=[a_{ij}]\in \mathbf{R}^{N\times N}$ is symmetric.
In $\mathcal{G}$, the $i$-th vertex represents the $i$-th agent. Let $a_{ij}>0$ if and only if $(i, j)\in \mathcal{E}$,
i.e., there is a communication link between agents $i$ and $j$. Undirected graph $\mathcal{G}$ is connected
if any two distinct agents of $\mathcal{G}$ can be connected via a path that follows the edges of $\mathcal{G}$.
For agent $i$, the degree is defined as $d_{i}\triangleq \sum_{j=1}^{N}a_{ij}$. Diagonal matrix
$\mathcal{D}\triangleq diag\{d_{1}, \ldots, d_{N}\}$ is used to denote the degree matrix of diagraph $\mathcal{G}$. Denote
the Laplacian matrix by $L_G=\mathcal{D}-\mathcal{A}$. The eigenvalues of $L_G$ are denoted by $\lambda_i\in R, i=1,\cdots, N,$ and an ascending order
in magnitude is written as $0=\lambda_1\leq \cdots\leq\lambda_N,$ that is, the Laplacian matrix $L_G$ of a undirected graph has at
least one zero eigenvalue and all the nonzero eigenvalues are in the
open right half plane. Furthermore, $L_G$ has exactly one zero eigenvalue
if and only if $G$ is connected\cite{godsil}.

In the sequel, it is assumed that the communication topology is connected.


\subsection{Regular Optimal Control for Single-agent System}

In this part, we study the infinite-horizon optimization problem for the single-agent system where the linear system is given by
\begin{eqnarray}
\dot{x}(t)&=&\bar{A}x(t)+\bar{B}u(t),\label{r11}
\end{eqnarray}
and the cost function:
\begin{eqnarray}
J&=&\int_0^\infty\Big(x'(t)\bar{Q}x(t)+u'(t)\bar{R}u(t)\Big)dt,\label{r12}
\end{eqnarray}
where $\bar{A},\bar{B},\bar{Q},\bar{R}$ are constant matrices with compatible dimensions. Moreover,
$\bar{Q}$ and $\bar{R}$ are nonzero semi-positive definite.

Let $\bar{P}$ satisfy the following algebraic Riccati equation
\begin{eqnarray}
0&=&\bar{A}'\bar{P}+\bar{P}\bar{A}+\bar{Q}-\bar{P}\bar{B}\bar{R}^{\dag}\bar{B}'\bar{P},\label{r5}\\
\bar{B}'\bar{P}&=&\bar{R}\bar{R}^{\dag}\bar{B}'\bar{P}.\label{r6}
\end{eqnarray}
With (\ref{r6}), the Riccati equation (\ref{r5}) is called regular and the corresponding LQ problem is named by
regular LQ problem \cite{Zhang18}. In this paper, we mainly consider the regular LQ problem.
By using the solvability of the Riccati equation (\ref{r5})-(\ref{r6}) and the stabilization of the specific matrix, we have the solvability of the infinite-horizon optimization problem.
\begin{theorem}\label{t2}
The optimal and stabilizing solution to minimize the cost function (\ref{r12}) exists if and only if the following two items hold:
\begin{enumerate}
  \item ARE (\ref{r5}) and (\ref{r6}) admits a solution $P\geq 0;$
  \item System $(\bar{A}-\bar{B}\bar{R}^{\dag}\bar{B}'\bar{P},\bar{B}(I-\bar{R}^{\dag}\bar{R}))$ is stabilizable.
\end{enumerate}
In this case, the optimal and stabilizing solution is given by
\begin{eqnarray}
u(t)&=&-\bar{R}^{\dag}\bar{B}'\bar{P}x(t)+(I-\bar{R}^{\dag}\bar{R})z(t),\label{r15}\\
z(t)&=&Kx(t),\label{r16}
\end{eqnarray}
where $K$ is chosen such that the matrix
\begin{eqnarray}
\bar{A}-\bar{B}\bar{R}^{\dag}\bar{B}'\bar{P}+\bar{B}(I-\bar{R}^{\dag}\bar{R})K\label{r19}
\end{eqnarray} is stable.
\end{theorem}
\emph{Proof.}
``Sufficiency" By taking derivative to $x'(t)\bar{P}x(t)$, we have
\begin{eqnarray}
\frac{d[x'(t)\bar{P}x(t)]}{dt}&=&x'(t)(\bar{A}'\bar{P}+\bar{P}\bar{A})x(t)+u'(t)\bar{B}'\bar{P}x(t)\nonumber\\
&&+x'(t)\bar{P}\bar{B}u(t).\nonumber
\end{eqnarray}
Thus, it is obtained that
\begin{eqnarray}
&&x'(T)\bar{P}x(T)-x'(0)\bar{P}x(0)\nonumber\\
&=&\int_0^T\Big(-x'(t)\bar{Q}x(t)-u'(t)\bar{R}u(t)\nonumber\\
&&+[u(t)+\bar{R}^{\dag}\bar{B}'\bar{P}x(t)]'\bar{R}[u(t)+\bar{R}^{\dag}\bar{B}'\bar{P}x(t)]\Big)dt.\nonumber
\end{eqnarray}
The cost function (\ref{r12}) can be reformulated as
\begin{eqnarray}
J&=&\lim_{T\rightarrow\infty}\int_0^T\Big(x'(t)\bar{Q}x(t)+u'(t)\bar{R}u(t)\Big)dt\nonumber\\
&=&x'(0)\bar{P}x(0)-\lim_{T\rightarrow\infty}x'(T)\bar{P}x(T)\nonumber\\
&&+\int_0^\infty[u(t)+\bar{R}^{\dag}\bar{B}'\bar{P}x(t)]'\bar{R}[u(t)+\bar{R}^{\dag}\bar{B}'\bar{P}x(t)]dt.\nonumber\\\label{r18}
\end{eqnarray}
Since $\bar{R}\geq 0,$  it is known that the optimal controller is given by
\begin{eqnarray}
0&=&u(t)+\bar{R}^{\dag}\bar{B}'\bar{P}x(t).\nonumber
\end{eqnarray}
Combining with Lemma \ref{lem}, it is obtained that the controller in (\ref{r15}) satisfies the above equation.
Together with (\ref{r16}), it holds that
\begin{eqnarray}
\lim_{t\rightarrow\infty}x(t)=0.\nonumber
\end{eqnarray}
Accordingly, the optimal cost is derived from (\ref{r18}) that
\begin{eqnarray}
J^*=x'(0)\bar{P}x(0).\nonumber
\end{eqnarray}

``Necessity" Consider the Riccati differential equation:
\begin{eqnarray}
0&=&\dot{\bar{P}}_T(t)+\bar{A}'\bar{P}_T(t)+\bar{P}_T(t)\bar{A}+\bar{Q}\nonumber\\
&&-\bar{P}_T(t)\bar{B}\bar{R}^{\dag} \bar{B}'\bar{P}_T(t),\label{r21}\\
\bar{B}'\bar{P}_T(t)&=&\bar{R}\bar{R}^{\dag}\bar{B}'\bar{P}_T(t).\label{r22}
\end{eqnarray}
with $\bar{P}(T)=0.$
By the solvability of the regular LQ problem \cite{Zhang18}, that is, the finite-horizon optimal control problem
by minimizing the cost function $J_T=\int_0^T\Big(x'(t)\bar{Q}x(t)+u'(t)\bar{R}u(t)\Big)dt$ subject to (\ref{r11})
is solvable, it follows that (\ref{r21}) and (\ref{r22}) admit a solution $\bar{P}_T\geq0.$ Moreover, the optimal cost is
\begin{eqnarray}
J^*_T=x'(0)\bar{P}_T(0)x(0).\label{r23}
\end{eqnarray}

Noting that $\bar{Q}\geq0, \bar{R}\geq0,$ we have
\begin{eqnarray}
J_{T_1}\leq J_{T_2}, ~\mbox{for}~T_1\leq T_2.\nonumber
\end{eqnarray}
This further implies that
\begin{eqnarray}
\bar{P}_{T_1}(0)\leq \bar{P}_{T_2}(0).\nonumber
\end{eqnarray}
Since the optimal and stabilizing controller exists, then the system is stabilizable which leads to that
there exists a stabilizing controller such that $\|x(t)\|^2+\|u(t)\|^2\leq e^{-\alpha t}x(0).$
Thus there exists a positive constant $c$ such that
\begin{eqnarray}
J_T\leq c.\nonumber
\end{eqnarray}
Together with (\ref{r23}), we have the uniform boundedness of the matrix $\bar{P}_T(t).$ Thus
there exists $\bar{P}$ such that
\begin{eqnarray}
\lim_{T\rightarrow\infty} \bar{P}_T(0)=\bar{P}.\nonumber
\end{eqnarray}
Noting the time-invariance of $\bar{P}_T(t),$ i.e.,
\begin{eqnarray}
\bar{P}_{T}(t)=\bar{P}_{T-t}(0),\nonumber
\end{eqnarray}
there holds that
\begin{eqnarray}
\lim_{T\rightarrow\infty} \bar{P}_T(t)=\lim_{T\rightarrow\infty} \bar{P}_{T-t}(0)=\bar{P}.\nonumber
\end{eqnarray}
Moreover, $\bar{P}$ obeys the equation (\ref{r5}) and (\ref{r6}).

Using (\ref{r18}) and $\bar{R}\geq0,$ we have the optimal controller is given by (\ref{r15}).
Accordingly, the system (\ref{r11}) becomes
\begin{eqnarray}
\dot{x}(t)&=&\Big(\bar{A}-\bar{B}\bar{R}^{\dag}\bar{B}'\bar{P}\Big)x(t)+\bar{B}(I-\bar{R}^{\dag}\bar{R})z(t).\nonumber
\end{eqnarray}
If $(\bar{A}-\bar{B}\bar{R}^{\dag}\bar{B}'\bar{P},\bar{B}(I-\bar{R}^{\dag}\bar{R}))$ is unstabilizable, then for any $z(t),$ the corresponding
state diverges. This is a contradiction to the existence of the optimal and stabilizing controller. The proof is now completed. \hfill $\blacksquare$

\begin{remark}
In the literature, there exist some sufficient conditions to guarantee the existence of the solution to (\ref{r5}).
In fact, noting that $\bar{R}\geq 0,$ there exists an orthogonal matrix $\bar{T}$ such that
\begin{eqnarray}
\bar{T}'\bar{R}\bar{T}=\left[
       \begin{array}{cc}
         \bar{R}_1 & 0 \\
         0 & 0 \\
       \end{array}
     \right]~~~ \mbox{and}~~~ \bar{B}\bar{T}=\left[
                               \begin{array}{cc}
                                 \bar{B}_1 & \bar{B}_2 \\
                               \end{array}
                             \right].\nonumber
\end{eqnarray}
Then (\ref{r5})-(\ref{r6}) is reformulated as
\begin{eqnarray}
0&=&\bar{A}'\bar{P}+\bar{P}\bar{A}+\bar{Q}-\bar{P}\bar{B}_1\bar{R}_1^{-1}\bar{B}_1'\bar{P},\label{r13}\\
\bar{B}_2'\bar{P}&=&0.\label{r14}
\end{eqnarray}
The equation (\ref{r13}) is a standard algebraic Riccati equation (ARE). By using the result in \cite{zhang1},
under the assumption that $(\bar{A},\bar{B}_1)$ is stabilizable, then there exists a positive semidefinite solution $\bar{P}\geq0$ to (\ref{r13}).
If $(\bar{A},\bar{B}_1)$ is stabilizable and $(\bar{A},\bar{Q})$ is observable, then there exists
a unique positive definite solution $\bar{P}>0$ for ARE (\ref{r5}) and the solution is stabilizing, i.e., the
matrix $\bar{A}-\bar{B}_1\bar{R}_1^{-1}\bar{B}_1'\bar{P}$ is stable.
\end{remark}

\section{Optimal Consensus Control}

We now consider the distributed optimization of multi-agent system (\ref{r1}).

Denote
$X(t)=\left[
       \begin{array}{ccc}
         x_1'(t) & \cdots & x_n'(t) \\
       \end{array}
     \right]'$ and $U(t)=\left[
       \begin{array}{ccc}
         u_1'(t) & \cdots & u_n'(t) \\
       \end{array}
     \right]',$
then the multi-agent system (\ref{r1}) is reformulated as
\begin{eqnarray}
\dot{X}(t)=\bar{A}X(t)+\bar{B}U(t),\label{r3}
\end{eqnarray}
where $\bar{A}=I_{N}\otimes A,\bar{B}=I_{N}\otimes B.$

The cost function $\sum_{i=1}^NJ_i$ can be reformulated as
\begin{eqnarray}
J&=&\int_0^\infty\Big(X'(t)\bar{Q}X(t)+U'(t)\bar{R}U(t)\Big)dt,\label{r4}
\end{eqnarray}
where $\bar{Q}=I_{N}\otimes Q,\bar{R}=I_{N}\otimes R.$

Define the following Riccati equation:
\begin{eqnarray}
0&=&A'P+PA+Q-PBR^{\dag}B'P,\label{r24}\\
B'P&=&RR^{\dag}B'P,\label{r25}
\end{eqnarray}
and let $\mathcal{A}=A-BR^{\dag}B'P$ and $\mathcal{B}=B(I-R^{\dag}R).$
By applying an appropriate coordinate transformation, there exists an invertible matrix $T_1$ such that
$T_1^{-1}\mathcal{A}T_1=diag\{\mathcal{A}_s, \mathcal{A}_u\}$ and $T_1^{-1}\mathcal{B}=\left[
                                                                      \begin{array}{c}
                                                                        \mathcal{B}_s \\
                                                                        \mathcal{B}_u \\
                                                                      \end{array}
                                                                    \right],
$ where $A_s$ is stable and $A_u$ is unstable, that is, all the eigenvalues of $\mathcal{A}_u$ are either on or outside the unit disk.
From \cite{zhang1}, we have the following solvability result.
\begin{lemma}\label{lem1}
Assume that $(\mathcal{A}_u,\mathcal{B}_u)$ is stabilizable,
then there exists a unique solution $\mathcal{P}_u>0$ to the following Riccati inequality
\begin{eqnarray}
0&=&\mathcal{A}_u'\mathcal{P}_u+\mathcal{P}_u\mathcal{A}_u-\mathcal{P}_u\mathcal{B}_u\mathcal{B}_u'\mathcal{P}_u+I.\nonumber
\end{eqnarray}
\end{lemma}
We are now in the position to state the main result of the distributed optimization control.
\begin{theorem}\label{t1}
The optimal and consensus solution to minimize $\sum_{i=1}^NJ_i$ subject to (\ref{r1}) exists if and only if
the following items hold:
\begin{enumerate}
  \item ARE (\ref{r24}) and (\ref{r25}) admits a solution $P\geq 0;$
  \item System $(\mathcal{A},\mathcal{B})$ is stabilizable.
\end{enumerate}
In this case, the optimal solution is given by
\begin{eqnarray}
u_i(t)&=&-R^{\dag}B'Px_i(t)+(I-R^{\dag}R)z_i(t),\label{r7}\\
z_i(t)&=&\mathcal{K}\sum_{j=1}^{N}a_{ij}\Big(x_j(t)-x_i(t)\Big),\label{r26}
\end{eqnarray}
where \begin{eqnarray}
\mathcal{K}&=&\left[
                  \begin{array}{cc}
                    0 & \mathcal{K}_u \\
                  \end{array}
                \right]T_1^{-1},\nonumber\\
\mathcal{K}_u&=&\max\{1, \frac{1}{\min_{i=2,\cdots,N}\lambda_i}\}\mathcal{B}_u'\mathcal{P}_u,\label{r20}
\end{eqnarray}
while $\lambda_i,i=2,\cdots,N$ is the eigenvalues of the Laplacian matrix $L_G.$
\end{theorem}
\emph{Proof.}
``Sufficiency"
By applying similar procedures to (\ref{r18}), it yields that
\begin{eqnarray}
J&=&\sum_{i=1}^{N}x_i'(0)Px_i(0)-\sum_{i=1}^N\lim_{T\rightarrow\infty}x_i'(T)Px_i(T)\nonumber\\
&&+\sum_{i=1}^{N}\int_0^\infty[u_i(t)+R^{\dag}B'Px_i(t)]'R\nonumber\\
&&\times [u_i(t)+R^{\dag}B'Px_i(t)]dt.\label{r28}
\end{eqnarray}
Since $R\geq 0,$ it is known that the optimal controller is given by
\begin{eqnarray}
0=u_i(t)+R^{\dag}B'Px_i(t).\nonumber
\end{eqnarray}
Combining with Lemma \ref{lem}, the controller given by (\ref{r7}) satisfies the above equation.
By substituting (\ref{r7}) into (\ref{r1}) yields that
\begin{eqnarray}
\dot{x}_i(t)&=&Ax_i(t)-BR^{\dag}B'Px_i(t)+B(I-R^{\dag}R)z_i(t)\nonumber\\
&=&\mathcal{A}x_i(t)+\mathcal{B}z_i(t),\label{r8}
\end{eqnarray}
We now substitute (\ref{r26}) with $\mathcal{K}$ defined in (\ref{r20}) into (\ref{r8}) and
denote $\delta(t)=\left[
                                                    \begin{array}{ccc}
                                                      \delta_2(t) & \cdots & \delta_N(t) \\
                                                    \end{array}
                                                  \right]
$ with $\delta_i(t)=x_i(t)-x_1(t),$ we have the following dynamic:
\begin{eqnarray}
\dot{\delta}(t)&=&(I_N\otimes\mathcal{A})\delta(t)-(L_{22}+\textbf{1}_{N-1}\alpha')\otimes (\mathcal{B}K)\delta(t),\nonumber
\end{eqnarray}
where $L_{22}=\left[
                \begin{array}{cccc}
                  d_2 & -a_{23} & \cdots & -a_{2N} \\
                  \cdots & \cdots & \cdots & \cdots \\
                  -a_{N2} & -a_{N3} & \cdots & d_{N} \\
                \end{array}
              \right]
$ and $\alpha=\left[
                \begin{array}{ccc}
                  a_{12} & \cdots & a_{1N} \\
                \end{array}
              \right]'.$
From \cite{ma}, there exists an invertible matrix $T_2$ such that $T_2^{-1}(L_{22}+\textbf{1}_{N-1}\alpha')T_2=diag\{\lambda_2,\cdots,\lambda_N\}.$
Let $\tilde{\delta}(t)=\left[
                                                    \begin{array}{ccc}
                                                      \tilde{\delta}_2(t) & \cdots & \tilde{\delta}_N(t) \\
                                                    \end{array}
                                                  \right]\triangleq T_2^{-1}\delta(t),$ we have
\begin{eqnarray}
\dot{\tilde{\delta}}(t)&=&\Big(I_N\otimes\mathcal{A}-diag\{\lambda_2, \cdots,\lambda_N\}\otimes (\mathcal{B}K)\Big)\delta(t),\nonumber
\end{eqnarray}
that is,
\begin{eqnarray}
\dot{\tilde{\delta}}_i(t)&=&\Big(\mathcal{A}-\lambda_i\mathcal{B}\mathcal{K}\Big)\tilde{\delta}_i(t),i=2,\cdots,N.\label{r10}
\end{eqnarray}
Further denote $\bar{\delta}_i(t)=T_1^{-1}\tilde{\delta}_i(t),$ the dynamic can be given from (\ref{r10}) by
\begin{eqnarray}
\dot{\bar{\delta}}_i(t)=\Big(\left[
                          \begin{array}{cc}
                            \mathcal{A}_s & 0 \\
                            0 & \mathcal{A}_u \\
                          \end{array}
                        \right]-\lambda_i\left[
                                           \begin{array}{c}
                                             \mathcal{B}_s \\
                                             \mathcal{B}_u \\
                                           \end{array}
                                         \right]\mathcal{K}T_1\Big)\bar{\delta}_i(t).\nonumber
\end{eqnarray}
Together with $\mathcal{K}=\left[
                  \begin{array}{cc}
                    0 & \mathcal{K}_u \\
                  \end{array}
                \right]T_1^{-1},
$ the above equation is reformulated as
\begin{eqnarray}
\dot{\bar{\delta}}_i(t)&=&\left[
                              \begin{array}{cc}
                                \mathcal{A}_s & -\lambda_i\mathcal{B}_s\mathcal{K}_u \\
                                0 & A_u-\lambda_i\mathcal{B}_u\mathcal{K}_u \\
                              \end{array}
                            \right]\bar{\delta}_i(t).\nonumber
\end{eqnarray}
Since the matrices $ A_s$ and $ A_u-\lambda_i\mathcal{B}_u\mathcal{K}_u $ are all stable, it yields that
\begin{eqnarray}
\lim_{t\rightarrow\infty} \bar{\delta}_i(t)=0.\nonumber
\end{eqnarray}
This gives that the multi-agent system (\ref{r1}) achieves consensus. In addition, let $y(t)=\frac{1}{N}\sum_{i=1}^Nx_i(t),$
it is obtained that
\begin{eqnarray}
\dot{y}(t)&=&\mathcal{A}y(t),\label{r27}
\end{eqnarray}
with initial value $y(0)=\frac{1}{N}\sum_{i=1}^Nx_i(0).$
Accordingly, the consensus value of (\ref{r1}) is given by
\begin{eqnarray}
\lim_{t\rightarrow\infty}\|x_i(t)-y(t)\|^2=0.\label{r30}
\end{eqnarray}

Denote
\begin{eqnarray}
V(t)=y'(t)Py(t),\nonumber
\end{eqnarray}
which is positive semi-definite.
Together with (\ref{r24}), it follows that
\begin{eqnarray}
\dot{V}(t)&=&y'(t)\Big(\mathcal{A}'P+P\mathcal{A}\Big)y(t)\nonumber\\
&=&-y'(t)(Q+PBR^{\dag}B'P)y(t)\nonumber\\
&\leq &0.\nonumber
\end{eqnarray}
Thus,
\begin{eqnarray}
\lim_{t\rightarrow\infty} V(t)=\lim_{t\rightarrow\infty} y'(t)Py(t)\triangleq \alpha,\nonumber
\end{eqnarray}
where $\alpha$ is a nonnegative constant.
Note that
\begin{eqnarray}
&&\Big\|x_i'(t)Px_i(t)-y'(t)Py(t)\Big\|^2\nonumber\\
&=&\Big\|\Big(x_i(t)-y(t)\Big)'P\Big(x_i(t)-y(t)\Big)\nonumber\\
&&+2y'(t)P\Big(x_i(t)-y(t)\Big)\Big\|^2\nonumber\\
&\leq&2\Big\|\Big(x_i(t)-y(t)\Big)'P\Big(x_i(t)-y(t)\Big)\Big\|^2\nonumber\\
&&+4\Big\|y'(t)P\Big(x_i(t)-y(t)\Big)\Big\|^2\nonumber\\
&\leq&2\Big\|\Big(x_i(t)-y(t)\Big)'P\Big(x_i(t)-y(t)\Big)\Big\|^2\nonumber\\
&&+4\Big(y'(t)Py(t)\Big)\Big\|P^{\frac{1}{2}}\Big(x_i(t)-y(t)\Big)\Big\|^2.\nonumber
\end{eqnarray}
Recalling (\ref{r30}), it is further obtained that
\begin{eqnarray}
\lim_{t\rightarrow\infty} x_i'(t)Px_i(t) =\lim_{t\rightarrow\infty} y'(t)Py(t)=\alpha.\nonumber
\end{eqnarray}
In view of (\ref{r28}), the optimal cost is given by
\begin{eqnarray}
J^*&=&\sum_{i=1}^{N}x_i'(0)Px_i(0)-\alpha N.\label{r29}
\end{eqnarray}

``Necessity" The proof of the necessity is similar to that of Theorem \ref{t2}.
So we omit it.
The proof is now completed. \hfill $\blacksquare$

\begin{remark}
\begin{itemize}
  \item When $R=0,$ the distributed optimization problem is reduced to the standard consensus problem. Theorem \ref{t1} can be reformulated as:
The consensus control of system (\ref{r1}) exists if and only if system $(A,B)$ is stabilizable.
In fact, when $R=0,$ from (\ref{r7}) and (\ref{r26}),
we have
\begin{eqnarray}
u_i(t)&=&\mathcal{K}\sum_{j=1}^{N}a_{ij}\Big(x_j(t)-x_i(t)\Big),\nonumber
\end{eqnarray}
where $\mathcal{K}=\max\{1, \frac{1}{\min_{i=2,\cdots,N}\lambda_i}\}B'\hat{P}$
while $\hat{P}$ satisfies the following Riccati equation
\begin{eqnarray}
0=A'\hat{P}+\hat{P}A+I-\hat{P}BB'\hat{P}.\nonumber
\end{eqnarray}
This is consistent with Theorem 2 in \cite{ma}.
  \item When $R$ is invertible, there seldom exists a optimal consensus controller in general since the optimal controller exists
  uniquely which may not possess a distributed structure. The only one to ensure the consensus is that the optimal controller guarantee the
  stability of all agents. In this case, the consensus value is zero.
\end{itemize}

\end{remark}

\section{Numerical Examples}

\textbf{Example 1:}
Consider the optimization problem where the matrices in (\ref{r1}) and (\ref{r2}) are given as
\begin{eqnarray}
A=\left[
    \begin{array}{cc}
      0 & 1 \\
      0 & 0 \\
    \end{array}
  \right],B=\left[
    \begin{array}{cc}
      0 & 0 \\
      0 & 1 \\
    \end{array}
  \right],\nonumber\\
Q=\left[
    \begin{array}{cc}
      0 & 0 \\
      0 & 0 \\
    \end{array}
  \right], R=\left[
    \begin{array}{cc}
      1 & 0 \\
      0 & 0 \\
    \end{array}
  \right].\nonumber
\end{eqnarray}
The solution to (\ref{r24}) and (\ref{r25}) is $P=\left[
    \begin{array}{cc}
      0 & 0 \\
      0 & 0 \\
    \end{array}
  \right].$ Thus,
\begin{eqnarray}
\mathcal{A}=\left[
    \begin{array}{cc}
      0 & 1 \\
      0 & 0 \\
    \end{array}
  \right], \mathcal{B}=\left[
    \begin{array}{cc}
      0 & 0 \\
      0 & 1 \\
    \end{array}
  \right].\nonumber
\end{eqnarray}
It is easily verified that $rank(\mathcal{B},\mathcal{A}\mathcal{B})=2,$
which implies that $(\mathcal{A},\mathcal{B})$ is controllable and thus stabilizable.
Let the communication topology of the multi-agent is given by Figure 1. The Lapalacian matrix is
as $\left[
      \begin{array}{cccc}
        2 & -1 & -1 & 0 \\
        -1 & 2 & -1 & 0 \\
        -1 & -1 & 3 & -1 \\
        0 & 0 & -1 & 1 \\
      \end{array}
    \right]
$ and the eigenvalues are in fact $0,1,3,4.$
\begin{center}
\begin{figure}\scalebox{0.45}{
\hs{0mm}\includegraphics{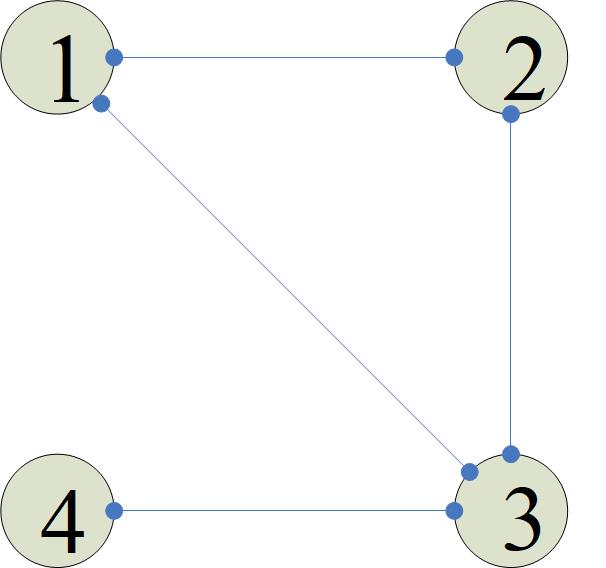}}
\caption{The communication topology graph.}
\end{figure}
\end{center}
By applying Theorem \ref{t1}, the optimal controller is
\begin{eqnarray}
u_i(t)&=&(I-R^{\dag}R)z_i(t)\nonumber\\
&=&(I-R^{\dag}R)\max\{1, \frac{1}{\min_{i=2,\cdots,N}\lambda_i}\}\mathcal{B}'\mathcal{P}\nonumber\\
&&\times \sum_{j=1}^{N}a_{ij}\Big(x_j(t)-x_i(t)\Big),\nonumber
\end{eqnarray}
where $\mathcal{P}$ is the solution to the following Riccati equation:
\begin{eqnarray}
0&=&\mathcal{A}'\mathcal{P}+\mathcal{P}\mathcal{A}-\mathcal{P}\mathcal{B}\mathcal{B}'\mathcal{P}+I.\nonumber
\end{eqnarray}
By solving the above Riccati equation, the optimal consensus protocol is in fact as
\begin{eqnarray}
u_i(t)&=&\max\{1, \frac{1}{\min_{i=2,\cdots,N}\lambda_i}\}\left[
    \begin{array}{cc}
      0 & 0 \\
      0 & 1 \\
    \end{array}
  \right]\left[
    \begin{array}{cc}
      1 & \sqrt{3} \\
      \sqrt{3} & 3 \\
    \end{array}
  \right]\nonumber\\
  &&\times \sum_{j=1}^{N}a_{ij}\Big(x_j(t)-x_i(t)\Big)\nonumber\\
&=&\left[
    \begin{array}{cc}
      0 & 0 \\
      \sqrt{3} & 3 \\
    \end{array}
  \right] \sum_{j=1}^{N}a_{ij}\Big(x_j(t)-x_i(t)\Big).
  \end{eqnarray}
From (\ref{r29}), we have that the optimal cost is zero. As shown in Fig. 2, the errors among
the states tend to zero. This gives the consensus of the multi-agent system.

\begin{center}
\begin{figure}\scalebox{0.45}{
\hs{0mm}\includegraphics{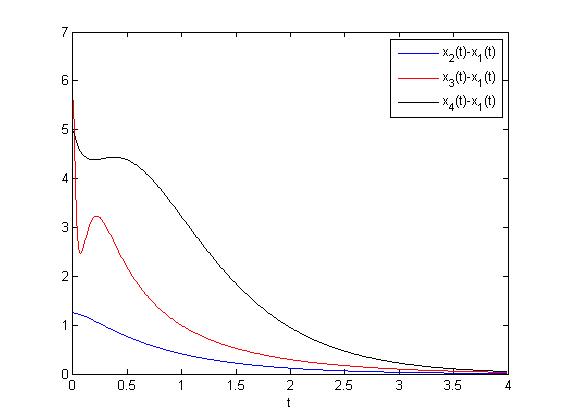}}
\caption{The trajectories of the errors among the states.}
\end{figure}
\end{center}

\textbf{Example 2:}
Consider the optimization problem where the matrices in (\ref{r1}) and (\ref{r2}) are given as
\begin{eqnarray}
A=0, B=\left[
         \begin{array}{cc}
           1 & 1 \\
         \end{array}
       \right],Q=1,R=\left[
    \begin{array}{cc}
      1 & 0 \\
      0 & 0 \\
    \end{array}
  \right].\nonumber
\end{eqnarray}
The solution to (\ref{r24}) and (\ref{r25}) is $P=1.$ Thus,
\begin{eqnarray}
\mathcal{A}=-1, \mathcal{B}=\left[
                              \begin{array}{cc}
                                0 & -1 \\
                              \end{array}
                            \right]
.\nonumber
\end{eqnarray}
It is easily verified that $rank(\mathcal{B},\mathcal{A}\mathcal{B})=1,$
which implies that $(\mathcal{A},\mathcal{B})$ is controllable and thus stabilizable.

Let the communication topology of the multi-agent is also given by Figure 1.
By applying Theorem \ref{t1}, the optimal controller is
\begin{eqnarray}
u_i(t)&=&\left[
            \begin{array}{c}
              -1 \\
              0 \\
            \end{array}
          \right]x_i(t)\nonumber.
\end{eqnarray}

In this case, the closed-loop system becomes
\begin{eqnarray}
\dot{x}_i(t)=-x_i(t),\nonumber
\end{eqnarray}
which implies that the multi-agent system achieves consensus and the consensus value is zero.
The optimal cost is accordingly $\sum_{i=1}^Nx_i'(0)x_i(0).$

\section{Conclusions}

In this paper, we studied the consensus control of multi-agent system with optimization performance. By considering
a cost function with a positive semi-definite weighting matrix of control, we derived the distributed and optimal control
protocol to ensure the consensus and optimization of the multi-agent system. In particular, the derived controller
is composed of two parts: One part is the feedback of the individual state which minimizes the cost function and the
other part is the feedback of the relative states between the neighbours which guarantees the consensus.


\end{document}